\documentclass[12pt]{article}
\def\demo{\noindent{\bf Proof .-}}
\newtheorem{theorem}{Theorem}
\newtheorem{lemma}{Lemma}
\newtheorem{proposition}{Proposition}
\newtheorem{definition}{Definition}
\newtheorem{example}{Example}
\newtheorem{remark}{Remark}

\begin{document}
\begin{center}
{\LARGE\bf {A note on Veronese varieties}}
\end{center}
\vskip.5truecm
\begin{center}
{Margherita Barile\footnote{Partially supported by PRIN Algebra Commutativa e Computazionale, Italian Ministry of Education,University and Research.}\\ Dipartimento di Matematica, Universit\`{a} di Bari, Via E. Orabona 4,\\70125 Bari, Italy\\e-mail:barile@dm.uniba.it}
\end{center}
\vskip1truecm
\noindent
{\bf Abstract} We show that for every prime $p$, there is a class of Veronese varieties which are set-theoretic complete intersections if and only if the ground field has characteristic $p$.\vskip.5truecm
\section*{Introduction} An affine or projective variety $V$ is called a {\it set-theoretic complete intersection} if it can be defined by the least possible number of equations, i.e., by codim\,$V$ equations. This property can depend upon the characteristic of the ground field, but not many examples of this kind have been discovered so far. In \cite{BL} the authors describe a class of codimension 2 toric varieties which are set-theoretic complete intersections only in one positive characteristic. In this paper we show that there is a class of Veronese varieties, of arbitrarily high codimension, having the same property. We shall present them by means of their parametrizations as toric varieties, which will enable us to apply the criteria on the associated semigroup developed in \cite{BMT}. We shall also use cohomological methods to provide lower bounds for the minimum number of defining equations.

\section{Preliminary results}
Let $K$ be an algebraically closed field. Let $p$ be a prime number, $h$ a positive integer and $n\geq3$ an integer.  Consider the following subset of ${\bf N}^n$:
$$T=\{\sum_{i=1}^na_i{\bf e}_i\vert\sum_{i=1}^n a_i=p^h,a_i\in{\bf N}\mbox{ for all }i=1,\dots, n\},$$
\noindent
where $\{{\bf e}_1,\dots, {\bf e}_n\}$ is the standard basis of ${\bf Z}^n$.   With every $\sum_{i=1}^na_i{\bf e}_i\in T$ we can associate the $p^h$-uple 
$$(i_1,\dots, i_{p^h})=(\underbrace{1,\dots, 1}_{a_1\ {\rm times }},\underbrace{2,\dots, 2}_{a_2\ {\rm times }},\dots, \underbrace{n,\dots, n}_{a_n\ {\rm times }}),$$
\noindent
and this defines a bijection between $T$ and the set
\begin{equation}\label{string}P=\{(i_1,\dots, i_{p^h})\mid 1\leq i_1\leq\cdots\leq i_{p^h}\leq n\}.\end{equation}
\noindent We know that $\vert T\vert={{n+p^h-1}\choose{p^h}}$.  In the affine space $K^{\vert T\vert}$ we fix the coordinates $x_{i_1\dots i_{p^h}}$, $(i_1,\dots, i_{p^h})\in P$. 
With $T$ we can associate the simplicial affine toric variety $V\subset K^{\vert T\vert}$,  defined by the following parametrization:
$$V=V^n_{p^h}:\left\{\begin{array}{rcl} x_{1\dots 1}&=&u_1^{p^h}\\
x_{2\dots 2}&=&u_2^{p^h}\\
&\vdots&\\
x_{n\dots n}&=&u_n^{p^h}\\
x_{11\dots112}&=&u_1^{p^h-1}u_2\\
x_{11\dots122}&=&u_1^{p^h-2}u_2^2\\
&\vdots&\\
x_{12\dots222}&=&u_1u_2^{p^h-1}\\
x_{11\dots123}&=&u_1^{p^h-2}u_2u_3\\
&\vdots&\\
x_{n-1n\dots n}&=&u_{n-1}u_n^{p^h-1}
\end{array}\right.
$$
\noindent 
It has codimension $N=\vert T\vert -n$. In the sequel, for the sake of simplicity, we shall set $p^h=q$.
Note that $V$ is the affine cone over the projective toric variety  of ${\bf P}^{\vert T\vert -1}$ having the same parametrization; the latter is called {\it Veronese variety}, and was extensively studied by Gr\"obner \cite{G}.  \par\smallskip\noindent
We consider the polynomial ring $R=K[x_{i_1\dots i_{q}}\mid (i_1,\dots, i_{q})\in P]$, and, in view of (\ref{string}), we define the {\it content} of the indeterminate $x_{i_1\dots i_{q}}$ as the set (with repeated elements):
$$\left\{\underbrace{u_1,\dots, u_1}_{a_1\ {\rm times }},\underbrace{u_2,\dots, u_2}_{a_2\ {\rm times }},\dots, \underbrace{u_n,\dots, u_n}_{a_n\ {\rm times }}\right\}.$$
\noindent More generally, the content of a product of indeterminates will be the (disjoint) union of the contents of its factors. Let $I(V)$ be the defining ideal of $V$ in $R$. Then $I(V)$ is generated by binomials. Clearly, given two monomials $M,N\in R$, we have that $M-N\in I(V)$ if and only if $M$ and $N$ have the same content. We have just established the following 
\begin{proposition}\label{generators} The binomials in $I(V)$ are the (non zero) differences\par\smallskip\noindent
{\rm (*)}$\qquad F=$\newline
$x_{i_1\dots i_{q}}x_{i_{q+1}\dots i_{2q}}\cdots x_{i_{(s-1)q+1}\dots i_{sq}}-
x_{i_{\sigma(1)}\dots i_{\sigma(q)}}x_{i_{\sigma(q+1)}\dots i_{\sigma(2q)}}\cdots x_{i_{\sigma((s-1)q+1)}\dots i_{\sigma(sq)}}$,\par\bigskip\noindent
where $s$ is a positive integer, and $\sigma$ any element of the symmetric group $S_{sq}$.
\end{proposition}
\begin{remark}{\rm  We have that $F=0$ if and only if $\sigma\left(\left\{i_{(k-1)q+1},\dots, i_{kq}\right\}\right)=\left\{i_{(k-1)q+1},\dots, i_{kq}\right\}$ for some $k\in\{1,\dots, s\}$. This occurs if and only if $F$ has a monomial factor}.\end{remark}
\begin{example} {\rm For  $p=2$, $h=1$, and $n=3$  the variety $V$ admits the following parametrization 
$$V:\left\{\begin{array}{rcl}
x_{11}&=&u_1^2\\
x_{22}&=&u_2^2\\
x_{33}&=&u_3^2\\
x_{12}&=&u_1u_2\\
x_{13}&=&u_1u_3\\ 
x_{23}&=&u_2u_3
\end{array}\right.$$
\noindent
The ideal $I(V)$ is minimally generated by the following six quadratic binomials:
$$x_{12}^2-x_{11}x_{22},\ x_{13}^2-x_{11}x_{33},\ x_{23}^2-x_{22}x_{33},$$
$$x_{12}x_{33}-x_{13}x_{23},\ x_{13}x_{22}-x_{12}x_{23},\ x_{23}x_{11}-x_{12}x_{13}.$$}
\end{example}
\noindent
In general we have
\begin{theorem}\label{degree2} $I(V)$ is generated by its binomials of degree 2.
\end{theorem}
\demo Let $B$ the set of binomials of type (*) having degree 2.  Let $F$ be any non zero binomial of type (*), and suppose that $\sigma=\tau_1\cdots\tau_m$, where, for all $i=1,\dots, m$, $\tau_i$ is a 2-cycle. In view of Remark 1, we can skip all  $\tau_i$ which involve two indices belonging to the same indeterminate.  Up to a change of indices we may thus assume that $\tau_m=(1\,q+1)$. We prove that $F\in(B)$ by induction on $m\geq1$. First suppose that $m=1$. Then $\sigma=(1\, q+1)$, so that
\begin{eqnarray*} F&=&x_{i_1i_2\dots i_{q}}x_{i_{q+1}i_{q+2}\dots i_{2q}}x_{i_{2q+1}\dots i_{3q}}\cdots x_{i_{(s-1)q+1}\dots i_{sq}}\\
&&-x_{i_{q+1}i_2\dots i_{q}}x_{i_1i_{q+2}\dots i_{2q}}x_{i_{2q+1}\dots i_{3q}}\cdots x_{i_{(s-1)q+1}\dots i_{sq}}\\
&=&(x_{i_1i_2\dots i_{q}}x_{i_{q+1}i_{q+2}\dots i_{2q}}-x_{i_{q+1}i_2\dots i_{q}}x_{i_1i_{q+2}\dots i_{2q}})x_{i_{2q+1}\dots i_{3q}}\cdots x_{i_{(s-1)q+1}\dots i_{sq}},
\end{eqnarray*}
\noindent
where the term in brackets belongs to $B$. Hence $F\in(B)$. Now assume that $m>1$ and suppose the claim true for $m-1$. We have:
\begin{eqnarray*} 
F&=&x_{i_1i_2\dots i_{q}}x_{i_{q+1}i_{q+2}\dots i_{2q}}x_{i_{2q+1}\dots i_{3q}}\cdots x_{i_{(s-1)q+1}\dots i_{sq}}\\
&&-x_{i_{q+1}i_2\dots i_{q}}x_{i_1i_{q+2}\dots i_{2q}}x_{i_{2q+1}\dots i_{3q}}\cdots x_{i_{(s-1)q+1}\dots i_{sq}}\\
&&+x_{i_{q+1}i_2\dots i_{q}}x_{i_1i_{q+2}\dots i_{2q}}x_{i_{2q+1}\dots i_{3q}}\cdots x_{i_{(s-1)q+1}\dots i_{sq}}\\
&&-x_{i_{\sigma(1)}i_{\sigma(2)}\dots i_{\sigma(q)}}x_{i_{\sigma(q+1)}i_{\sigma(q+2)}\dots i_{\sigma(2q)}}x_{i_{\sigma(2q+1)}\dots i_{\sigma(3q)}}\cdots x_{i_{\sigma((s-1)q+1)}\dots i_{\sigma(sq)}}\\
&=&(x_{i_1i_2\dots i_{q}}x_{i_{q+1}i_{q+2}\dots i_{2q}}-x_{i_{q+1}i_2\dots i_{q}}x_{i_1i_{q+2}\dots i_{2q}})x_{i_{2q+1}\dots i_{3q}}\cdots x_{i_{(s-1)q+1}\dots i_{sq}}\\
&&+x_{i_{q+1}i_2\dots i_{q}}x_{i_1i_{q+2}\dots i_{2q}}x_{i_{2q+1}\dots i_{3q}}\cdots x_{i_{(s-1)q+1}\dots i_{sq}}\\
&&-x_{i_{\sigma\tau_m(q+1)}i_{\sigma\tau_m(2)}\dots i_{\sigma\tau_m(q)}}x_{i_{\sigma\tau_m(1)}i_{\sigma\tau_m(q+2)}\dots i_{\sigma\tau_m(2q)}}\\
&&\qquad\qquad\qquad\qquad\qquad x_{i_{\sigma\tau_m(2q+1)}\dots i_{\sigma\tau_m(3q)}}\cdots x_{i_{\sigma\tau_m((s-1)q+1)}\dots i_{\sigma\tau_m(sq)}}.
\end{eqnarray*}
\noindent
Induction applies to $\sigma\tau_m=\tau_1\cdots\tau_{m-1}$, so that $F\in(B)$ in this case, too. This completes the proof.
\par\smallskip\noindent
The previous lemma was already shown in \cite{G}, but our approach is simpler and emphasizes the combinatorial aspect. In \cite{G} one can also find a proof of the next result, for which we follow a more direct method. \par\smallskip\noindent 
Let $\{{\bf e}_{i_1\dots i_{q}}\mid (i_1,\dots, i_{q})\in P\}$ be the standard basis of ${\bf Z}^{\vert T\vert}.$ Here we assume that the indices $(i_1,\dots, i_{q})$ are arranged in the ascending lexicographic order. 
\begin{lemma}\label{singularity} The origin is the only singular point of $V$. 
\end{lemma}
\demo We prove the claim using the Jacobian criterion. Let 
\begin{eqnarray*}
\!\!\!\!\!{\bf w}&=&
(\bar x_{1\dots1},\dots,\bar x_{n\dots n},\bar x_{1\dots12},\dots,\bar x_{n-1n\dots n})\\
&=&(\bar u_1^{q},\dots, \bar u_n^{q},\bar u_1^{q-1}\bar u_2,\dots, \bar u_{n-1}\bar u_n^{q-1})\in V,
\end{eqnarray*}
 \noindent and let ${\bf J}({\bf w})$ be a Jacobian matrix associated with the set $B$ of generators, evaluated at ${\bf w}$. Its rows are
$$\sum_{(i_1,\dots, i_{q})\in P}\left(\frac{\partial F}{\partial x_{i_1\dots i_{q}}}\right)({\bf w}){\bf e}_{i_1\dots i_{q}},$$
\noindent where $F\in B$. Clearly ${\bf J}(0)$ is the zero matrix, so that $0$ is a singular point for $V$. Now suppose that ${\bf w}\ne0$. Up to a change of indices we may assume that $\bar u_1\ne0$. We prove that rank\,${\bf J}({\bf w})\geq N$. Set
$$F_{i_1\dots i_{q}}=x_{1\dots1}x_{i_1\dots i_{q}}-x_{1\dots1 i_{q}}x_{1i_1\dots i_{q-1}},$$
\noindent
for all 
$$(i_1,\dots,i_{q})\in P'=P\setminus\{(1,\dots,1,i)\mid i=1,\dots, n\}.$$
\noindent
These are $N=\vert T\vert-n$ elements of $B$. Let ${\bf J}'$ be the $N\times N$-submatrix of ${\bf J}({\bf w})$ formed by the entries
$$\left(\frac{\partial F_{i_1\dots i_{q}}}{\partial x_{j_1\dots j_{q}}}\right)({\bf w})\qquad\mbox{for all }(i_1,\dots, i_{q}), (j_i,\dots, j_{q})\in P'.$$
\noindent The rows of ${\bf J}'$ are
$$\bar x_{1\dots1}{\bf e}_{i_1\dots i_{q}}-\bar x_{1\dots1 i_{q}}{\bf e}_{1i_1\dots i_{q-1}},$$
\noindent where the second summand is missing whenever $i_1=\cdots= i_{q-2}=1$. Otherwise $(1,i_1,\dots, i_{q-1})$ is lexicographically smaller than $(i_1,\dots, i_{q})$. Hence ${\bf J}'$ is a lower triangular matrix whose diagonal entries are all equal to $\bar x_{1\dots1}\ne0$. Thus ${\bf J}'$ is invertible. It follows that rank\,${\bf J}({\bf w})\geq N$, and ${\bf w}$ is not a singular point for $V$. This completes the proof.
\par\smallskip\noindent
The following result contains the notion of Galois covering: we refer to \cite{M}, p.~43 for a definition.  
\begin{lemma}\label{Galois}
Suppose that char\,$K\ne p$. The map
$$\phi:K^n\setminus\{0\}\longrightarrow V\setminus\{0\}$$
$$\qquad\qquad\qquad\qquad(u_1,\dots, u_n)\longrightarrow (u_1^{q},\dots,u_n^{q},u_1^{q-1}u_2,\dots,u_{n-1}u_n^{q-1})$$
\noindent
is a finite Galois covering with cyclic Galois group of order $q$.
\end{lemma}
\demo We show that $\phi$ is an \'etale finite surjective map and that the  multiplicative group $C_{q}$ of the complex $q$-th roots of unity acts transitively and faithfully on every fiber of $\phi$.
Surjectivity and finiteness are clear. The coordinate ring of $V$ is  $A=K[u_1^{q},\dots,u_n^{q},u_1^{q-1}u_2,\dots,u_{n-1}u_n^{q-1}]$, and the coordinate ring of $K^n$ is
$$K[u_1,\dots, u_n]=A[T_1,\dots, T_n]/(T_1-u_1,\dots, T_n-u_n).$$
Set $P_i=T_i-u_i$ for all $i=1,\dots, n$. Then the Jacobian matrix $\left(\frac{\partial P_i}{\partial T_j}\right)$ is the identity matrix of order $n$. According to \cite{M}, Corollary 3.16, p.~27, this implies that $\phi$ is \'etale.  
\newline
Let ${\bf w}=(\bar u_1^{q},\dots,\bar u_n^{q},\bar u_1^{q-1}\bar u_2,\dots,\bar u_{n-1}\bar u_n^{q-1})\in V\setminus\{0\}$. Up to a permutation of indices we may assume that $\bar u_1=\cdots=\bar u_k=0$ and $\bar u_{k+1},\dots,\bar  u_n\ne0$, for some $k$, $0\leq k\leq n-1$. Let ${\bf v}=(0,\dots, 0, \bar u_{k+1},\dots, \bar u_n) \in\phi^{-1}({\bf w})$. Then for all ${\bf v'}\in \phi^{-1}({\bf w})$, ${\bf v'}=(0,\dots, 0, g_{k+1}\bar u_{k+1},\dots, g_n\bar u_n)$, where, for all $i=k+1,\dots, n$, $g_i\in C_{q}$. We are going to show that
\begin{equation}\label{g}  g_i=g_j\qquad\mbox{for all indices }i,j\mbox{ such that }k+1\leq i<j\leq n.\end{equation}
\noindent
Let $k+1\leq i<j\leq n$. First assume that $p=2$, and $h=1$. In this case $g_i\in\{-1,1\}$ for all $i=k+1,\dots, n$.
Equating the entries of $\phi({\bf v}')$ and $\phi({\bf v})$ of index 
$ij$ we have  $g_i\bar u_ig_j\bar u_j=\bar u_i\bar u_j$ from which we deduce that
$$g_ig_j=1,$$
\noindent and this implies (\ref{g}).
Now suppose that $p=2$ and $h>1$. 
Equating the entries of $\phi({\bf v}')$ and $\phi({\bf v})$ of index 
$$\underbrace{i\dots i}_{2^{h-1}}\underbrace{j\dots j}_{2^{h-1}}$$
\noindent we get $g_i^{2^{h-1}}\bar u_i^{2^{h-1}}g_j^{2^{h-1}}\bar u_j^{2^{h-1}}=\bar u_i^{2^{h-1}}\bar u_j^{2^{h-1}}$, whence
\begin{equation}\label{1}g_i^{2^{h-1}}g_j^{2^{h-1}}=1,\end{equation}
\noindent and equating their entries of index
$$\underbrace{i\dots i}_{2^{h-1}-1}\underbrace{j\dots j}_{2^{h-1}+1}$$
\noindent we get, similarly,
\begin{equation}\label{2}g_i^{2^{h-1}-1}g_j^{2^{h-1}+1}=1.\end{equation}
\noindent
If we divide (\ref{1}) by (\ref{2}) we obtain $g_ig_j^{-1}=1$, which implies (\ref{g}).
Finally suppose that $p$ is odd. Equating the entries of $\phi({\bf v}')$ and $\phi({\bf v})$ of index 
$$\underbrace{i\dots i}_{\frac{q+1}2}\underbrace{j\dots j}_{\frac{q-1}2}$$
\noindent we deduce that
\begin{equation}\label{3}g_i^{\frac{q+1}2}g_j^{\frac{q-1}2}=1,\end{equation}
\noindent and equating their entries of index
$$\underbrace{i\dots i}_{\frac{q-1}2}\underbrace{j\dots j}_{\frac{q+1}2}$$
\noindent we get
\begin{equation}\label{4}g_i^{\frac{q-1}2}g_j^{\frac{q+1}2}=1.\end{equation}
\noindent
From (\ref{3}) and (\ref{4}) we again deduce (\ref{g}).
It follows that, in all cases
$${\bf v}'=g{\bf v}\qquad\mbox{ for some }g\in C_{q}.$$
\noindent
This shows that 
$$\phi^{-1}({\bf v})\subset\{g{\bf v}\mid g\in C_{q}\},$$
\noindent
The opposite inclusion is obvious. This completes the proof.

 \section{The set-theoretic complete intersection property}
In this section we show that $V=V^n_q$ is a set-theoretic complete intersection (i.e., it is set-theoretically defined by $N$ equations) if and only if char\,$K=p$. We shall consider the subsemigroup ${\bf N}T$ of ${\bf N}^n$ and the subgroup ${\bf Z}T$ of ${\bf Z}^n$ generated by $T$. We need to recall the following  two definitions, both quoted from \cite{BMT}, pp.~1894--1895.
\begin{definition}\label{definition1}{\rm
Let $p$ be a prime number and let $T_1$ and $T_2$ be non-empty
subsets of $T$ such that $T     =   T_1\cup T_2$ and $T_1\cap T_2    =
\emptyset$.
Then $T$ is called a {\it p-gluing} of $T_1$ and
$T_2$ if  there is $s\in{\bf N}$ and  a
nonzero element ${\alpha}\in {{\bf N}}^n$ such that  ${{\bf Z}} T_1\cap{{\bf Z}} T_2    =  {{\bf Z}} {\alpha}$
and 
$  p^s{\alpha}\in {{\bf N} } T_1 \cap {{\bf N} } T_2$. \\}
\end{definition}
\begin{definition}\label{definition2}{\rm
An affine semigroup ${{\bf N} }T$ is called {\it completely $p$-glued} if $T$
is the
$p$-gluing of
$T_1$ and $T_2$, where each  of the semigroups ${{\bf N} } T_1, {{\bf N} }
T_2$
is completely $p$-glued or a free abelian semigroup.\\}
\end{definition}
There is a large class of completely $p$-glued semigroups:
\begin{lemma}\label{glued} Let $p$ be a prime number, $h$ a positive integer, and let
$$T_0=\{p^h{\bf e}_{1\dots1},\dots, p^h{\bf e}_{n\dots n}\}.$$
\noindent If $T_0\subset T\subset{\bf N}^n$, where $T$ is finite, then the semigroup ${\bf N}T$ is completely $p$-glued.
 \end{lemma}
\demo  We proceed by induction on $\vert T\vert=t\geq n$. Since $T_0$ is free, for $t=n$ there is nothing to prove. Suppose $t>n$ and the claim true for all smaller $t$. Pick any $\alpha\in T\setminus T_0$ and set $T_1=T\setminus\{\alpha\}$, $T_2=\{\alpha\}$. Then, being $T_0\subset T_1$, ${\bf N}T_1$ is completely $p$-glued by induction and, moreover, $p^h\alpha\in{\bf N}T_0\subset {\bf N}T_1$, so that $p^h\alpha\in{\bf N}T_1\cap {\bf N}T_2$. Hence $T$ is the $p$-gluing of $T_1$ ad $T_2$. This proves that ${\bf N}T$ is completely $p$-glued.
\par\bigskip\noindent
In the same paper (see \cite{BMT}, Theorem 5) we find the following characterization of set-theoretic complete intersections:
\begin{theorem}\label{theorem1}
Suppose that char\,$K=p>0$. Then a toric variety in over $K$ is a set-theoretic complete intersection on binomials if and only if the associated semigroup is completely $p$-glued. 
\end{theorem}
We shall also need the following criterion, cited from \cite{BS}, Lemma 3$^\prime$. 
\begin{lemma}\label{Newstead}Let
$W\subset\tilde W$ be affine varieties. Let $d=\dim\tilde
W\setminus W$. If there are $s$ equations $F_1,\dots, F_s$ such
that $W=\tilde W\cap V(F_1,\dots,F_s)$, then 
$$H^{d+i}_{\rm et}(\tilde W\setminus W,{\bf Z}/r{\bf Z})=0\quad\mbox{ for all
}i\geq s$$ and for all $r\in{\bf Z}$ which are prime to {\rm char}\,$K$.
\end{lemma}
Finally we recall one result from group cohomology. 
\begin{lemma}\label{group} Let $p$ be a prime, and $h$ be a positive integer.  Let $G$ be a cyclic group of order $p^h$. Then  
$$H^i(G,{\bf Z}/p^h{\bf Z})\ne0\quad\mbox{for all }i\geq0,$$
with respect to any $G$-action on ${\bf Z}/p^h{\bf Z}$.
\end{lemma} 
\demo We consider $G$ as a multiplicative group. Let $g$ be a generator of $G$, and let  $g[1]=[a]$,  where $a\in{\bf Z}$ and [] denotes the residue class mod $p^h$. Then the action of $g$ on $A={\bf Z}/p^h{\bf Z}$ is multiplication by $a$. Moreover, $g^{p^h}=1$, so that
$[1]=g^{p^h}[1]=[a]^{p^h}$. Hence  $a^{p^h}\equiv 1$ (mod $p$). By Fermat's Little Theorem it follows that $a\equiv 1$ (mod $p$). Consequently,  $p^{h-1}a\equiv p^{h-1}$ (mod $p^h$). Thus
$$g[p^{h-1}]=[p^{h-1}a]=[p^{h-1}].$$
This shows that $[p^{h-1}]\in A^G$, the submodule of $g$-invariants of $A$. Now, by definition of group cohomology, $H^0(G;A)=A^G$,  so that $H^0(G;A)\ne 0.$ By a well-known property (see, e.g., \cite{W}, Proposition 3-2-3), all cohomology groups of a cyclic group with respect to a finite $G$-module $A$ have the same order. This completes the proof. 
\par\bigskip\noindent
We are now ready to prove:
\begin{theorem}\label{theorem2}  $V$ is a set-theoretic complete intersection if and only if char\,$K=p$.
\end{theorem}
\demo Let $T$ be the subset of ${\bf N}^n$ associated with $V$. Then ${\bf N}T$ is completely $p$-glued by virtue of Lemma \ref{glued}.  By Theorem \ref{theorem1} it follows that, if char\,$K=p$, $V$ is a set-theoretic complete intersection. Now suppose that char\,$K\ne p$. According to \cite{M}, Theorem 2.20, p.~105, the map $\phi$ in
Lemma \ref{Galois} gives rise to the following Hochschild-Serre spectral sequence:
$$H^i({\bf Z}/q{\bf Z},H^j_{\rm et}(K^n\setminus\{0\},{\bf Z}/q{\bf Z}))\Longrightarrow H^{i+j}_{\rm et}(V\setminus\{0\},{\bf Z}/q{\bf Z}),$$
\noindent
where $H^{\cdot}$ and $H_{\rm et}^{\cdot}$ denote group cohomology and \'etale cohomology respectively. Recall that
\begin{equation}\label{point} H_{\rm et}^j(K^n\setminus\{0\},{\bf Z}/q{\bf Z})\simeq\cases{{\bf Z}/q{\bf Z}& if $j=0, 2n-1$\cr
0& else}\end{equation}
\noindent
and that, according to Lemma \ref{group},
\begin{equation}\label{cyclic} H^i({\bf Z}/q{\bf Z},{\bf Z}/q{\bf Z})\ne0\qquad\mbox{ for all }i\geq0.\end{equation}
\noindent
Since $2n-1>2$, (\ref{point}) and (\ref{cyclic}) 
imply that 
\begin{equation}\label{ne0}H_{\rm et}^{n+1}(V\setminus\{0\},{\bf Z}/q{\bf Z})\simeq H^{n+1}({\bf Z}/q{\bf Z},H_{\rm et}^0(K^n\setminus\{0\},{\bf Z}/q{\bf Z}))\ne0.\end{equation}
\noindent
Since, by Lemma \ref{singularity}, $V\setminus\{0\}$ is a smooth $n$-dimensional variety, we can apply Poincar\'e Duality to it (see \cite{M}, Theorem 1.11, p.~276), so that
\begin{equation}\label{poincare}H_c^{n-1}(V\setminus\{0\},{\bf Z}/q{\bf Z})\simeq H_{\rm et}^{n+1}(V\setminus\{0\},{\bf Z}/q{\bf Z}),\end{equation}
\noindent where $H_c^{\cdot}$ denotes cohomology with compact support. In the sequel we shall omit the coefficient group ${\bf Z}/q{\bf Z}$ for the sake of simplicity.
 From (\ref{ne0}) and (\ref{poincare}) we conclude that
\begin{equation}\label{6}H_c^{n-1}(V\setminus\{0\})\ne0.\end{equation}
\noindent 
We have a long exact sequence of cohomology with compact support
$$\cdots\longrightarrow H_c^{n-2}(0)\longrightarrow H_c^{n-1}(V\setminus\{0\})\longrightarrow H_c^{n-1}(V)\longrightarrow H_c^{n-1}(0)\longrightarrow\cdots,$$
\noindent
where, being $n\geq3$, $H_c^{n-2}(0)=H_c^{n-1}(0)=0$. Hence 
$$H_c^{n-1}(V)\simeq H_c^{n-1}(V\setminus\{0\}),$$
\noindent
 so that, by (\ref{6}),
\begin{equation}\label{7}H_c^{n-1}(V)\neq0.\end{equation}
We also have a long exact sequence
$$\cdots\longrightarrow H_c^{n-1}(K^{N+n})\longrightarrow H_c^{n-1}(V)\longrightarrow H_c^{n}(K^{N+n}\setminus V)\longrightarrow H_c^{n}(K^{N+n})\longrightarrow\cdots,$$
\noindent
where $H_c^{n-1}(K^{N+n})=H_c^{n}(K^{N+n})=0$, so that
$H_c^{n}(K^{N+n}\setminus V)\simeq H_c^{n-1}(V)$, and thus, in view of (\ref{7}), 
$$H_c^{n}(K^{N+n}\setminus V)\ne0.$$
\noindent 
If we apply Poincar\'e Duality, we finally obtain
$$H_{\rm et}^{N+n+N}(K^{N+n}\setminus V)\ne0.$$
By Lemma \ref{Newstead} it follows that, under our present assumption that char\,$K\ne p$, $V$ is not set-theoretically defined by $N$ equations, i.e. it is not a set-theoretic complete intersection. This completes the proof.
\par\bigskip\noindent
According to the methods developed in \cite{BMT},  $V^n_{q}$ is, in characteristic $p$, a set-theoretic complete intersection on the following set of $N$ binomials
$$x_{i_1\dots i_{q}}^{q}-x_{1\dots1}^{i_1}\cdots x_{q\dots q}^{i_{q}},\qquad(i_1,\dots, i_p)\in P\setminus\{(1\dots1),\dots,(n\dots n)\}.$$
\noindent The same set of $N$  binomials was found by Gattazzo \cite{Ga} by direct computations.
\begin{remark}{\rm 
The variety $V^{\frac{n(n+1)}2}_2$ (of codimension $\frac{n(n-1)}2$) is the determinantal variety defined by the vanishing of the 2-minors of an $n\times n$ symmetric variety of indeterminates over $K$. Theorem \ref{theorem1} generalizes part of the results in \cite{B}: there it was shown that the minimum number of equations required to define $V^{\frac{n(n+1)}2}_2$ set-theoretically is 
$$\cases{\frac{n(n+1)}2-2&if char\,$K\ne2$\cr
\frac{n(n-1)}2&if char\,$K=2$}
$$
\noindent
which means that, for $n\geq3$,  $V^{\frac{n(n+1)}2}_2$ is a set-theoretic complete intersection if and only if char\,$K=2$. }
\end{remark}

\begin{remark}{\rm A general lower bound for the minimum number of equations which define a variety set-theoretically (the so-called {\it arithmetical rank}, ara) is given by the  {\it local cohomological dimension}: if $I=I(V)$, this number is
$${\rm cd}\,I=\max\{n\in{\bf N}\mid H_I^n(R)\ne0\},$$
\noindent
where $H_I^{\cdot}$ denotes local cohomology with respect to $I$. In the case of the Veronese variety $V=V^n_{p^h}$, the ideal $I$ is perfect (see \cite{G}, p.~259), so that, according to \cite{PS}, Prop.~4.1, 
$${\rm cd}I={\rm ht}\,I=N.$$
\noindent However, it follows from Theorem \ref{theorem2} that
$$N={\rm cd}I<{\rm ara}\,I \qquad\mbox{ if char}\,K\ne p,$$
\noindent
 so that the lower bound for the arithmetical rank provided by the local cohomological dimension is almost always non-sharp. Varieties where the local cohomological dimension and the arithmetical rank differ are particularly sought. The only previously known examples are the Pfaffian ideals of an alternating matrix of indeterminates, in all positive characteristics (see \cite{B}, Remarks 6.2), and Reisner's variety, in all characteristics different from 2 (see  \cite{SV}, p.~250, \cite{L}, Example 1, and \cite{Y}, Example 2).}
\end{remark}


\begin{thebibliography}{BMT}
\bibitem{B} Barile, M.: Arithmetical ranks of ideals associated to symmetric and alternating matrices. J.~Algebra, \textbf{176}, 59--82 (1995)
\bibitem{BL} Barile, M., Lyubeznik, G.: Set-theoretic complete intersections in characteristic $p$. Preprint (2004). To appear in: Proceedings of the American Mathematical Society. 
\bibitem{BMT}  Barile, M., Morales, M., Thoma, A.: Set-Theoretic Complete Intersections on Binomials. Proc.~Amer.~Soc., \textbf{130}, 1893--1903 (2002) 
\bibitem{BS} Bruns, W., Schw\"anzl: The number of
equations defining a determinantal variety. Bull.~London Math.~Soc., \textbf{22}, 439--445 (1990)
\bibitem{Ga} Gattazzo, R.: In characteristic $p=2$ the Veronese variety $V^m\subset{\bf P}^{m(m+3)/2}$ and each of its generic projection  is set-theoretic complete intersection, in: Greco, S., Strano, R. (eds.), Complete Intersections, Acireale 1983, Lecture Notes in Mathematics \textbf{1092}, 221--228, Springer, Berlin-Heidelberg (1984)
\bibitem{G} Gr\"obner, W.: \"Uber Veronesesche Variet\"aten und deren Projektionen. Arch.~Math., \textbf{16}, 257--264 (1965)
\bibitem{L} Lyubeznik, G.: On the local cohomology modules $H^i_{\cal A}(R)$ for ideals ${\cal A}$ generated by monomials in an $R$-sequence, in: Greco, S., Strano, R. (eds.), Complete Intersections, Acireale 1983, Lecture Notes in Mathematics \textbf{1092}, 221--228, Springer, Berlin-Heidelberg (1984)
\bibitem{M} Milne, J.: \'Etale Cohomology. Princeton University Press, Princeton (1980)
\bibitem{PS} Peskine, C., Szpiro, L.: Dimension finie et cohomologie locale. Inst.~Hautes \'Etudes Sci.~Publ.~Math.~ \textbf{42}, 47--119 (1973) 
\bibitem{SV} Schmitt, Th., Vogel, W.: Note on Set-Theoretic Intersections of Subvarieties of Projective Space. Math.~Ann. \textbf{245}, 247--253 (1979)
\bibitem{W} Weiss, E.: Cohomology of Groups. Academic Press, New York, London (1969)
\bibitem{Y} Yan, Z.: An \'etale analog of the Goresky-MacPherson formula for subspace arrangements. J.~Pure Appl.~Algebra \textbf{146}, 305--318 (2000)


\end{thebibliography}
\end{document}